\documentclass[12pt]{article} %CTEX报告文章格式

\usepackage[top=3cm,bottom=2cm,left=2cm,right=2cm]{geometry} % 页边距
\usepackage{amsmath} %数学公式
\usepackage{amsthm}
\usepackage{amssymb}
\usepackage{graphicx} %图片
\usepackage{tikz} %画图
\usepackage{listings}
\usepackage{amsfonts}
\usepackage{subfigure}
\usepackage{float}
\usepackage[all]{xy}

\usepackage{float}%提供float浮动环境
\usepackage{booktabs}%提供命令\toprule、\midrule、\bottomrule
\usetikzlibrary{shapes,arrows} %tikz图形库
\usepackage{overpic} %图上标记
\usepackage{authblk}
\usepackage{multirow}
\theoremstyle{thmstyletwo}%

\theoremstyle{thmstylethree}%
\usepackage{graphics}
\usepackage{epstopdf}
\usepackage{makecell}
\begin{document}
%\CTEXoptions[contentsname={\bfseries\zihao{4} Contens}]
%\setcitestyle{square,aysep={},yysep={;}}
%\captionwidth{0.8\textwidth}
%\changecaptionwidth
%\thispagestyle{empty}
%\pagestyle{plain}
%\newpage
%\setcounter{page}{1}
%\pagenumbering{Roman}
%\noindent\addcontentsline{toc}{section}{ÕªÒª}

\title{Spiking hierarchy in an adaptive exponential integrate-and-fire network synchronization}

%\title[Spiking hierarchy]{Spiking hierarchy in an adaptive exponential integrate-and-fire network synchronization}

\author[1,+]{Xiaoyue Wu}
\author[1,2,3,+]{Congping Lin}
\author[1,2,3,*]{Yiwei Zhang}
\affil[1]{\small School of Mathematics and Statistics, Huazhong University of Sciences and
Technology, Wuhan, China}
\affil[2]{\small Center for
Mathematical Sciences, Huazhong University of Science and
Technology, Wuhan, China}
\affil[3] {\small Hubei Key Laboratory of Engineering Modeling
and Scientific Computing, Huazhong University of Science and
Technology, Wuhan, China}
\affil[+] {\small These authors contributed equally to this work}
\affil[*] {\small Corresponding author {\tt{yiweizhang@hust.edu.cn}}}

\date{\today}
%\tableofcontents
%\setcounter{page}{1}
%\pagenumbering{arabic}
%\pagestyle{headings}
\maketitle

%\author[1]{\fnm{Xiaoyue} \sur{Wu}}%\email{}%\email{m201970072@hust.edu.cn}
%\author[1,2,3]{\fnm{Congping} \sur{Lin}}\email{congpinglin@hust.edu.cn}
%\author*[1,2,3]{\fnm{Yiwei} \sur{Zhang}}\email{yiweizhang@hust.edu.cn}
%%\equalcont{These authors contributed equally to this work.}
%
%\affil*[1]{\orgdiv{School of Mathematics and Statistics}, \orgname{Huazhong University of Science and
%Technology}, \orgaddress{ \city{Wuhan}, \postcode{430074}, \state{Hubei}, \country{China}}}
%\affil*[2]{\orgdiv{Center for Mathematical Sciences}, \orgname{Huazhong University of Science and
%Technology}, \orgaddress{\city{Wuhan}, \postcode{430074}, \state{Hubei}, \country{China}}}
%\affil*[3]{\orgdiv{Hubei Key Laboratory of Engineering Modeling
%and Scientific Computing}, \orgname{Huazhong University of Science and
%Technology}, \orgaddress{\city{Wuhan}, \postcode{430074}, \state{Hubei}, \country{China}}}

%\tableofcontents
%\setcounter{page}{1}
%\pagenumbering{arabic}
%\pagestyle{headings}
\begin{abstract}
Neuronal network synchronization has received wide interests. Network connection structure is known to play a key role in its synchronization. In the present manuscript, we study the influence of initial membrane potentials together with network topology on bursting synchronization, in particular the sequential spiking order in stabilized inter bursts. We find a hierarchical phenomenon on spiking order. We grade neurons into different layers: primary neurons are those initiate the spike and lower-layer neurons are then determined via network connections successively. This constructs a directed graph to indicate spiking propagation among different layers of neurons. Neurons in upper layers spike earlier than those in lower layers. More interestingly, we find that among the same layer, neuron spiking order is mainly associated with stimuli they receive from upper layer via network connections; more stimuli leads to earlier spiking. Furthermore, we find that with effectively the same stimuli from upper layer, neurons with less connections spike earlier.
\end{abstract}

%\date{\today}
%%\keywords{Spiking hierarchy, neuronal network, synchronization, bursting pattern}
%\maketitle

\section{Introduction}
The study of brain rhythms and synchronization of oscillatory activity has attracted substantial attention. In the brain, there are abundant experimental evidences of neuronal synchronization, which is associated to cognitive processes, such as visual cognition \cite{CTB}, memory formation~\cite{AMFE} and directed attention~\cite{MDG}.

%Synchronization can occur between coupled chaotic systems \cite{PC}.
A number of mathematical models have been developed on the basis of some physiological experiments to describe neuronal activities. The Hudgkin-Huxley model was proposed to simulate the potential behavior of cell membranes, based on equivalent circuit principles and the results of experiments with the giant axons of gun squid \cite{HH}. Its simplified versions including, Leaky Integrate-and-Fire neuron models, Izhikevich neuron models, and the adaptive exponential integrate-and-fire neuron model (aEIF) have been widely investigated \cite{Izh,LIF,TKK,BG}. In particular, the aEIF model has an exponential spiking mechanism combined with an adaptation \cite{BG}. Protachevicz {\em et al.} verified that there exist spike and burst synchronization in different cortical areas in a random network of aEIF neurons \cite{PBRB}. Based on the aEIF model, Borges {\em et al.} verified that bursting synchronization was more robust than spiking synchronization \cite{BPL}. More recently, in random networks of aEIF neural system, cluster synchronization induced by time delay has also been studied \cite{LYZ}.

Neuronal network synchronization can be influenced by dynamical properties of individual neurons. Naud {\em et al.} have shown multiple firing patterns of the aEIF model and studied different firing types in association with individual neuron dynamics \cite{NMC}. Two typical rhythm synchronization of coupled neurons: spiking synchronization which is equivalent to phase synchronization and bursting synchronization which occurs more often than spiking synchronization, are studied in \cite{SL}.

Besides individual neuron dynamics, network connections are also shown to play a key role in neuronal synchronization. To simulate the neuronal dynamic behavior in complex brain network, different network models have been proposed, including small-world network \cite{WS2,WS3,WS4}, Erd\"{o}s-R\'{e}nyi network \cite{ER}, scale-free network \cite{ScaleFree1,ScaleFree2}. In the Erd\"{o}s-R\'{e}nyi network, it was suggested that large connection probability promotes bursting synchronization \cite{BPL}. In small-world networks, adding some random connections instead of increasing rewiring probability, would increase neuronal synchronization~\cite{CGF}. The state of synchronization can be various including complete synchrony, delay synchrony and phase locking \cite{WLC}. Considering the dynamical heterogeneity, one may have different synchronization properties, even for networks with the same propensity to synchronization \cite{LSB}. In pulse-coupled networks of bursting Hindmarsh-Rose neurons, the stability of completely synchronous state in such networks was found to solely depend on the number of signals each neuron received \cite{BdH}.

We are interested in firing patterns of regular burst, in particular the firing sequence of neurons in phase locking. We find spiking hierarchy among neurons in inter-burst and study its associations with initial membrane potentials and neuronal network connection structures. We grade neurons in different layers, and find among the same layer, neuron receiving more stimuli from the upper-layer would spike earlier; receiving effectively the same stimuli from the upper-layer, neurons with less connections would spike earlier.

The manuscript is organised as follows: Sec.~\ref{sec:network model} describes the aEIF neuronal network model and shows parameter region for bursting and spiking synchronization. Network synchronization is investigated in Sec.~\ref{subsec:rgn} for regular networks and in Sec.~\ref{subsec:swn} for small-world networks. Conclusion are given in Sec.~\ref{sec:conclusion}.

\section{Neuronal network model}\label{sec:network model}
In this section, we first introduce the aEIF model and show parameter space giving burst and spike synchronization patterns in Sec.~\ref{sec:model} and then introduce the network model in Sec.~\ref{subsec:burst}.

\subsection{The aEIF model}\label{sec:model}
The aEIF model is a two-dimensional spiking neuron model given by the following differential equations \cite{BG}:
\begin{eqnarray}\label{eq:aEIF}
%\nonumber
%\Biggl\{
\left\{
\begin{aligned}
C\frac{dV}{dt} &=-g_{L}(V-E_{L})+g_{L}\triangle_{T}\exp(\frac{V-V_{T}}{\triangle_{T}})-w+I,\\
\tau_{w}\frac{dw}{dt} &=a(V-E_{L})-w,
\end{aligned}
\right.
\end{eqnarray}
where parameters are related to physiological quantities: $V(t)$ is the membrane potential, $w$ is the adaptation current, $I$ is the injected current, $C$ is the membrane capacitance, $g_{L}$ is the leak conductance, $E_{L}$ is the resting potential, $V_{T}$ is the threshold potential, $\triangle_{T}$ is the slope factor, $\tau_{w}$ is the time constant and $a$ is the level of threshold adaptation. Due to the exponential term, when the membrane potential $V$ is high enough, the trajectory will quickly diverge. When $V(t)$ reaches the threshold $V_T$, the membrane potential $V$ is reset to $V_{r}$ (referred as reset potential) and the adaptation current $w$ is increased by $b$ as in \cite{TB}:
\begin{align*}
V &\rightarrow V_{r},\\
w &\rightarrow w+b.
\end{align*}
%\end{equation}

In this aEIF model, many different firing patterns could occur including tonic spiking, regular bursting and adaptation \cite{BPL,Touboul}. In this manuscript we focus on regular bursting (i.e. the potential converges to a stable spiking cycle containing a specific number of spikes); a delay will occur before spiking sequence occurs regularly \cite{Touboul}. It is known that the injected current $I$ and the reset potential $V_r$ play important roles in generating different firing patterns; the magnitude of $I$ will affect the position of $V$-nullcline and $w$-nullcline of Eq.~(\ref{eq:aEIF})~\cite{BPL,Touboul}. We mainly consider influence of reset potential $V_{r}$ and injected current $I$ when fixing other parameters $C=281.0~pF, g_{L}=30.0~nS, E_{L}=-70.6~mV, V_{T}=-50.4~mV, \triangle_{T}=2~mV, \tau_{w}=20~ms, a=4~nA, b=0.5~nA$, as set in the \cite{TB}.

\begin{figure}[h]
	%\centering{\includegraphics[scale=0.13]{hebing1.pdf}}
    \centering
    \includegraphics[width=\columnwidth]{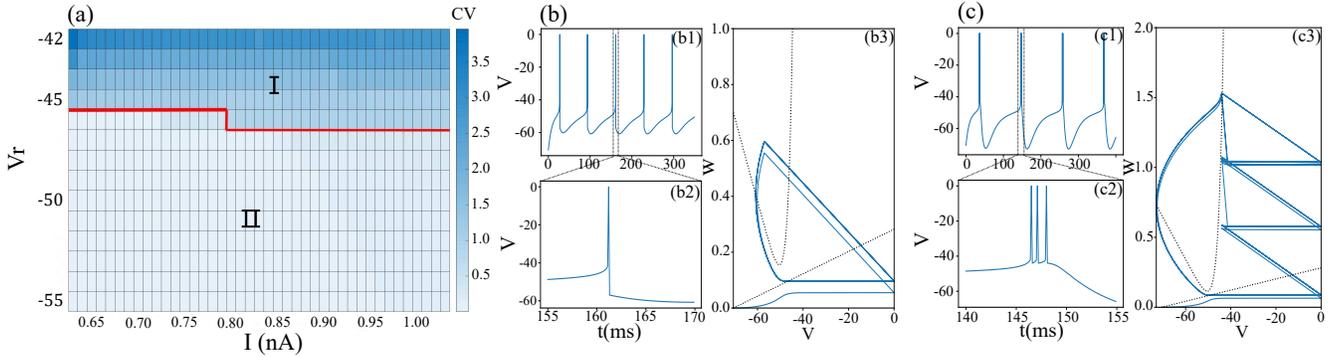}
	\caption{(a) shows the parameter space for two firing patterns as a function of the reset potential $V_{r}$ and input current $I$. Bursting pattern (region I, $CV\geq 0.5$) and spiking pattern (region II, $CV<0.5$) are separated by the red line. (b-c) shows examples of spiking and bursting pattern respectively; (b1-c1) show the time evolution of membrane potential (V) with zoom in (b2-c2) for a small time window; (b3-c3) show corresponding trajectories in $V$-$w$ phase; the trajectory in (b3) converges to a limit spiking cycle containing only one spike while the trajectory in (c3) converges to a stable spiking cycle containing 3 spikes. Dashed lines in (b3, c3) represent the $V$-nullcline and $w$-nullcline of Eq.~(\ref{eq:aEIF}). For visual display of curves in (b,c) we set $V_{T}=0~mV$ and $I=0.66nA$, $V_r=-44mV$ in (b), $I=0.70nA$, $V_r=-57mV$ in (c).}
	\label{fig:CV}
\end{figure}

We employ coefficient of variation ($CV$) of neuronal inter-spikes interval ($ISI$) \cite{GK} to distinguish different firing patterns:
\begin{equation}\label{eq:CV}
CV=\frac{\sigma_{ISI}}{\overline{ISI}},
\end{equation}
where $\sigma_{ISI}$ is the standard deviation of ISI normalized by its mean $\overline{ISI}$. If $CV<0.5$, the neuron fires in a spiking pattern; while $CV\geq0.5$, it fires in a bursting pattern~\cite{BPL}. Fig.~\ref{fig:CV}(a) shows that two types of firing occurs in $V_r$-I parameter space according to $CV$ given in Eq.~(\ref{eq:CV}). When the reset membrane potential $V_r$ is small, spiking patterns occur, whereas when $V_r$ is large, bursting patterns occur \cite{TB,Touboul}. Fig.~\ref{fig:CV}(b-c) show examples of tonic spiking and regular bursting patterns; in the burst pattern, it fires a specific number (3 here) of spikes per burst. With an input current, the neuron keeps periodic burst pattern after a certain period of adjustment as seen from trajectories shown in Fig.~\ref{fig:CV}(c3).

\subsection{aEIF neuronal network model and synchronization}\label{subsec:burst}
Next, we consider the aEIF neuronal network coupled through chemical synapses via the following differential equations~\cite{BPL}:
\begin{equation}
\left\{
\begin{aligned}
C\frac{dV_{i}}{dt}=&-g_{L}(V_{i}-E_{L})+g_{L}\triangle_{T}\exp(\frac{V_{i}-V_{T}}{\triangle_{T}})-w_{i}+I\\
&+(V_{REV}-V_{i})\sum_{j=1}^{N}M_{ij}g_{j},\\
\tau_{w}\frac{dw_{i}}{dt}=&a(V_{i}-E_{L})-w_{i},\\
\tau_{g}\frac{dg_{i}}{dt}=&-g_{i},
\end{aligned}
\right.
\label{eq:network}
\end{equation}
where $V_{i}, w_i, g_i$ represent the membrane potential, adaptation current and synaptic conductance of neuron $i$, $V_{REV}$ is the synaptic reversal potential, $\tau_{g}$ is the synaptic time constant, $M$ is the adjacency matrix with elements $M_{ij}=1$ or $0$ depending on whether or not there exists coupling between neurons $i$ and $j$. The synaptic conductance decays exponential with a synaptic time constant $\tau_g$. We consider connectivity between neurons is given by excitatory synapses. When $V_{i}$ reaches the threshold, the state variables are updated as follows~\cite{PBLL}:
%\begin{equation}
%\nonumber
\begin{align*}
V_{i} &\rightarrow V_{r},\\
w_{i} &\rightarrow w_{i}+b,\\
g_{i} &\rightarrow g_{i}+g_{exc},
\end{align*}
%\end{equation}
where $g_{exc}=0.05~nS$ as set in \cite{PBLL}. Eq.~(\ref{eq:network}) indicates that in addition to the external input current $I$, synaptic current is also generated by interactions between neurons. Thus in a network of $N$ neurons, the total current affected to neuron $i$ is:
\begin{align*}
I^{i}_{total}(t)=-w_{i}+I+(V_{REV}-V_{i})\sum_{j=1}^{N}M_{ij}g_{j}.
\end{align*}
We take the choice of single neuron aEIF model parameters in Section \ref{sec:model}, fix $\tau_{g}=2.728~ms$, $V_{REV}=0~mV$ as in \cite{PBLL}, and set initial values of $w$ and $g$ as $0$. Moreover, we fix $I=0.66~nA$ which will consecutively generate action potential, and $V_r=-44~mV$ where firing patterns with 3 spikes per burst is robust for input current $I$ in Fig.~\ref{fig:CV}(a). We study the effect of initial $V_i$'s on network synchronization. Note that initial values of $V_i$ are chosen between resting potential and threshold potential (i.e. $V_i\in [E_{L},V_{T}]$).

\subsubsection*{Synchronization parameter}
It is crucial to introduce parameters to identify neuronal synchronization. We introduce a synchronization parameter $\delta$ which focuses on the bursting patterns with specific number ($K$) of spikes per burst. We record the time of each spike (denoted as $t_{i,m-1}$ for the $m$-th spike of neuron $i$) and characterize synchronization degree at $n$-th burst (with $K$ spikes) as:
%\begin{equation}%\label{eq:delta_n}
%\nonumber
\begin{align*}
\delta(n):=\frac{1}{K}\sum_{j=0}^{K-1}{\frac{\sqrt{\frac{1}{N}{\sum_{i=0}^{N-1}(t_{i,j+K(n-1)})^2}-(\frac{1}{N}{\sum_{i=0}^{N-1}t_{i,j+K(n-1)}})^2}}{N-1}},
%\end{equation}
\end{align*}
where $N$ is the number of neurons in the network. This synchronization parameter $\delta(n)$ averages standard deviations of spiking time over all $K$ spikes in each burst and is adapted from \cite{ZL}. For parameters we focus, the system generates regular burst with 3 spikes per burst, i.e. $K=3$. Moreover, when the synchronization parameter $\delta(n)$ stabilize we define its stabilized synchronization parameter of network as:
\begin{align*}%
\label{eq:delta}
\nonumber
\delta:=\lim_{n\to+\infty}\delta(n).
\end{align*}
Smaller $\delta$ indicates spiking time are more clustered and synchronized.

\subsubsection*{Spiking time difference}
%It can be seen from the above discussion that when the synchronization parameter $\delta(n)$ is stable, the firing period of neurons in the coupled network will stabilize at the same value, as seen in Fig.~\ref{fig:DT}(a). For single neuron, when the parameters except the initial membrane potentials are given, the firing period of the neuron will be stable after a number of bursts regardless of its initial potential, as seen in Fig.~\ref{fig:DT}(b).

\begin{figure}[H]
	\centering
	\includegraphics[width=0.5\columnwidth]{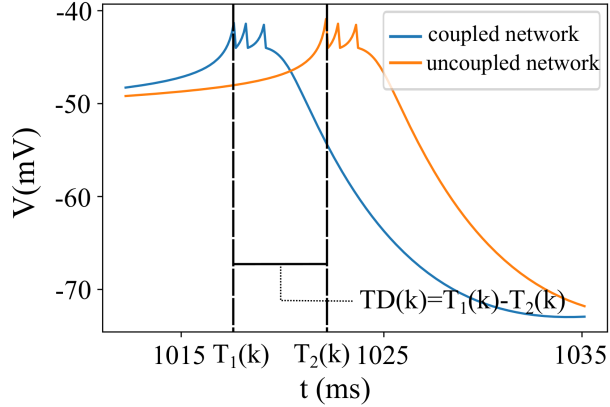}
	\caption{Illustration of time series of a neuron in the coupled (blue) or uncoupled network (orange) with certain $k$-th spiking time $T_{1,2}(k)$ and a time difference $TD(k)=T_1(k)-T_2(k)$.}
	\label{fig:DT}
\end{figure}

In order to quantify the effect of network coupling on spiking rhythm, for each $k$-th spike and each neuron, we calculate its spiking time $T_1(k)$ in the coupled network and also the spiking time $T_2(k)$ in the same scenario except that the network is uncoupled. We then calculate the spiking time difference $TD(k)=T_1(k)-T_2(k)$ between coupled and uncoupled networks as sketched in Fig.~\ref{fig:DT}. A spiking time difference $TD(k)<0$ (resp. $>0$) indicates that the network connections advance (resp. suppress) response of its action potential, while $TD(k)=0$ indicates that network connections have effectively no effect on its spiking time. Once the firing pattern and corresponding inter-burst period of neurons in the network reach stabilized, $TD(k)$ remains constant and we denote this stabilized time difference as $TD:=\lim_{k\to+\infty}TD(k)$. We use $TD_i(k)$ and $TD_i$ for the (stabilized) spiking time difference of particular neuron $i$ at $k$-th spike.

%\section{Explore the influencing factors of network synchronization effect}\label{sec:syn}
\section{Results}\label{subsec:ini_vol}
We first investigate in Sec.~\ref{subsec:ini_delta} \& \ref{subsec:rgn} the influence of initial membrane potentials on network synchronization using regular network where each neuron is connected to equal number of neighbors. With the assistance of spiking time difference calculation, we find spiking hierarchy of neurons in inter burst. We then examine this spiking hierarchy phenomenon in small world network in Sec.~\ref{subsec:swn} and discuss the spiking sequence order in relation with degree and initial membrane potentials of neurons within the same layer in the hierarch structure.

\subsection{Influence of initial potential on synchronization}\label{subsec:ini_delta}
\begin{figure}[H]
%\centering\includegraphics[scale=0.21]{deltaall11.pdf}%{p2delta.png}
\centering\includegraphics[width=0.4\columnwidth]{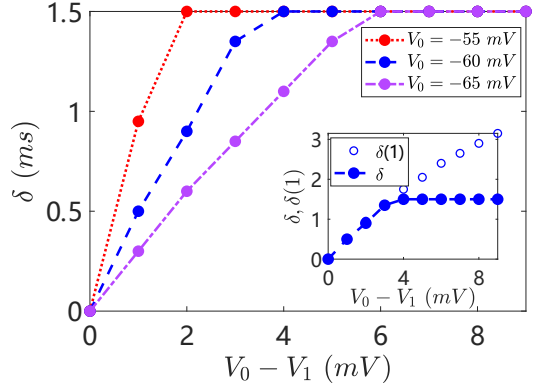}%{p2delta.png}
\caption{Stabilized synchronization parameter $\delta$ in a toy model with two nodes for different initial membrane potential $V_0$ and different initial potential difference $V_0-V_1$. The insert panel shows both $\delta(1)$ and stabilized synchronization parameter $\delta$ for $V_0=-60~mV$.}
\label{fig:toy}
\end{figure}

%\subsubsection*{A toy example: two-node network}
To investigate the interaction between neurons in the coupled neuronal network, we start with a toy model with two coupled nodes of initial membrane potentials $V_0$ and $V_1$ respectively. Without loss of generality, we consider $V_0\geq V_1$. Fig.~\ref{fig:toy} shows the stabilized synchronization parameter $\delta$ changes with the initial potential difference $V_0-V_1$. Clearly, $\delta=0$ when $V_0=V_1$ as expected. When $V_0-V_1$ increases, $\delta$ increases to a saturate value $1.5\times 10^{-3}~s$. The saturated synchronization parameter is independent of $V_{0,1}$ as seen from Fig.~\ref{fig:toy}.

We observe from Fig.~\ref{fig:toy} that for fixed initial potential difference $V_0-V_1$, larger $V_0$ gives higher $\delta$. Moreover, when increasing the potential difference $V_0-V_1$, larger $V_0$ can reach the saturated value with smaller $V_0-V_1$.

We show both synchronization parameter at first burst $\delta(1)$ and stabilized synchronization parameter $\delta$ in the insert figure in Fig.~\ref{fig:toy} and observe that $\delta(1)\approx\delta$ for small $V_0-V_1$. With large initial potential difference $V_0-V_1$ (where saturated synchronization parameter is reached), stabilized synchronization parameter $\delta$ could be much smaller than $\delta(1)$.

%with the increase of the initial membrane potential difference between the two neurons ($VD$), the initial synchronization become lager and the stable synchronization gradually approach a stable value. Note that, for different $V_0$, the value of stable synchronization are the same. Obviously, there are different critical value for $V_0$. If the initial potentials difference between two neurons exceeds this critical value, the stable synchronization parameters are the same regardless of the size of the $VD$. If VD is less than the critical value, there will be two situation: one is that the two neurons keep their initial firing state unchanged, so their synchronization also remains unchanged; The other is that the neurons in the network will adjust their firing state quickly and then reach a stable state, such as case in Fig.~\ref{fig:3delta}(a) with $V_0=-60mV$ and $VD=2mV$.

\begin{figure}[H]
%\centering\includegraphics[scale=0.21]{deltaall11.pdf}%{p2delta.png}
\centering\includegraphics[width=0.8\columnwidth]{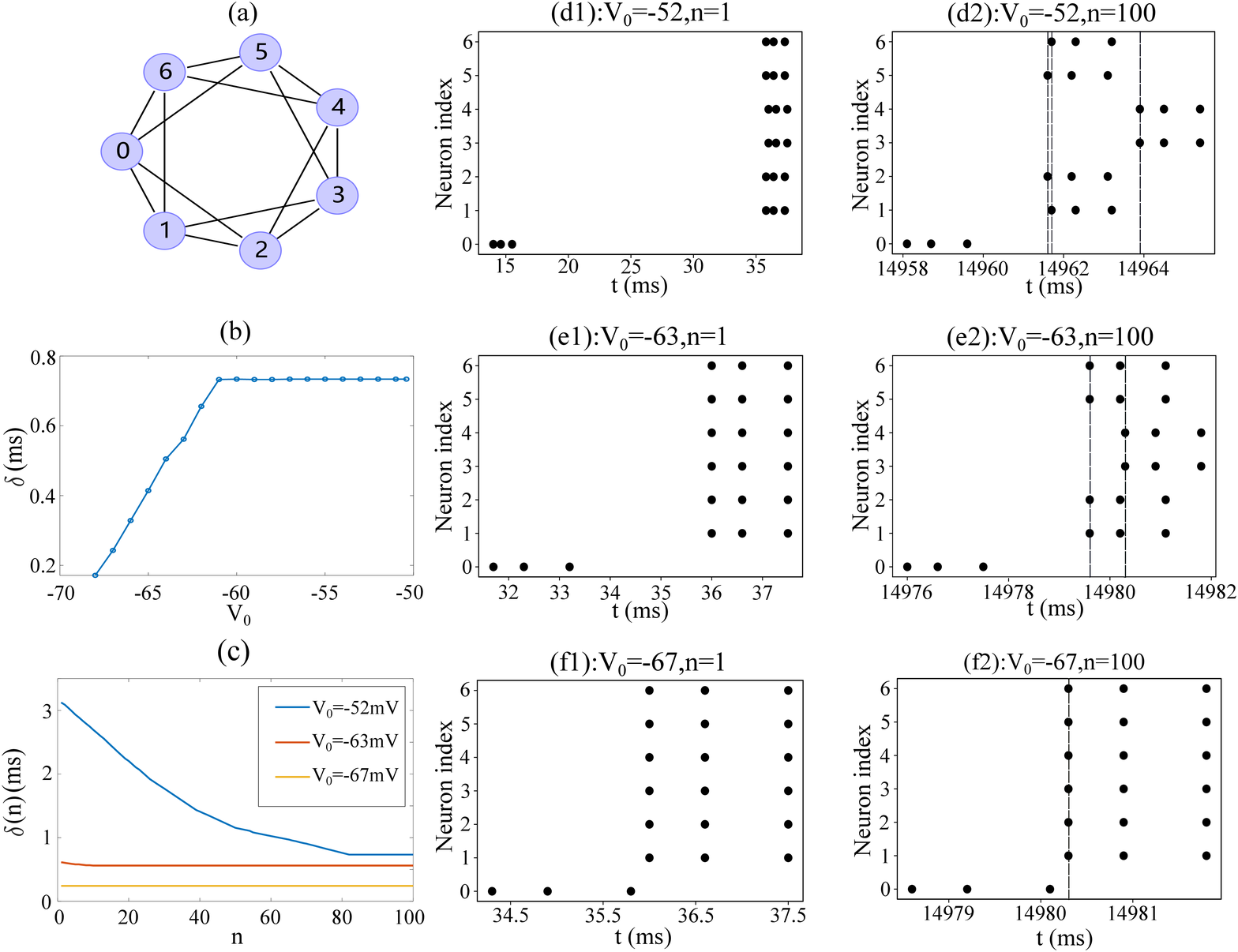}%{p2delta.png}
\caption{(a) shows a regular network of 7 neurons where each neuron is connected to its nearest 4 neighbors. (b) shows stabilized synchronization parameter $\delta$ when varying initial membrane potential $V_0$ while other potentials $V_i=-70.6~mV$ the same as resting potential in the regular network (a). (c) shows the evolution of synchronization parameter $\delta(n)$ for three initial membrane potential $V_0$, with the corresponding raster plots of excitatory neurons at $n=1$ and $n=100$ shown in panels (d2-f2).}
\label{fig:3delta}
\end{figure}

%\subsubsection*{regular network}
Next, we investigate the initial potential effect on neuronal synchronization in regular network with more nodes. We start with a simple case where initial membrane potentials are the same except for one neuron, say neuron 0 with initial potential $V_0$. Fig.~\ref{fig:3delta}(b) shows the stabilized synchronization parameter $\delta$ for various $V_0$ while other initial potentials are $V_i=-70.6~mV, i\neq 0$ (the same as resting potential $E_L$) in a  regular network of 7 neurons shown in Fig.~\ref{fig:3delta}(a). In this example, when $V_0<-61~mV$, the stabilized synchronization parameter $\delta$ exhibit approximately linear increase with $V_0$. When $V_0$ increases, the initial potential difference $V_0-V_i$ increases as well. While for larger $V_0$ ($V_0>-61mV$), $\delta$ saturates and keeps unchanged when increasing $V_0$. This is consistent with saturation observed in the toy example in Fig.~\ref{fig:toy}.

%When spiking pattern is stabilized, the state of excitatory neurons reach phase locking, i.e., firing sequence of excitatory neurons keeps the same.

Fig.~\ref{fig:3delta}(c) shows the evolution of synchronization parameter $\delta(n)$ for three different $V_0$, and the corresponding locked phases at two different bursts ($n=1$ and $n=100$) can be seen from raster plots shown in Fig.~\ref{fig:3delta}(d2-f2). From Fig.~\ref{fig:3delta}(b), we note that synchronization parameter $\delta$ stabilizes faster for smaller $V_0$. This is probably due to that the initial membrane potential of neuron 0 (which spike first among all neurons in the inter-burst) is closer to other initial potentials. Transition from the first burst to stabilized burst for these three different $V_0$ can be seen from Fig.~\ref{fig:3delta}(d-f). In all three different $V_0$, neuron 0 spike first in the inter-bursts, which then promote other neurons to spike .%Moreover, from the locked phase shown in Fig.~\ref{fig:3delta}(b2-d2), we observe that neuron $0$ always generate action potentials first in each burst, this might due to that $V_0$ is the largest initial membrane potential among all neurons.

Moreover, we observed in Fig.~\ref{fig:3delta}(d2,e2) that the spiking time in inter-bursts are different even for neurons with the same initial membrane potentials, due to neuronal network connections. In this example, neuron 0 has the largest initial potential and generates the action potential first. Then the neuron 0 affects its neighboring neurons $\{1,2,5,6\}$ via coupling and leads to different spiking time between neurons $\{1,2,5,6\}$ and $\{3,4\}$ as seen in Fig.~\ref{fig:3delta}(d2,e2). For $V_0=-67~mV$, we calculate $TD_i=0$ for all $i$, meaning neurons spike as effectively independent neurons, and those with larger initial potentials spike earlier than others. For $V_0=-63~mV$, $TD_{0,3,4}=0$ whereas $TD_{1,2,5,6}<0$ which is consistent with neurons $\{1,2,5,6\}$ spike earlier than neurons $\{3,4\}$, though they have the same initial potentials. For $V_0=-52~mV$, $TD_0=0$ and $TD_i<0$ for all $i\neq 0$, meaning that spiking time of neurons are advanced expect for neuron 0 due to network connection. It is worth noting that when $V_0>-61mV$ all states of phase locking are exactly the same.

\subsection{Spiking hierarchy in regular network synchronization}\label{subsec:rgn}
\begin{figure}[H]
%\centering
\centering
\includegraphics[width=\columnwidth]{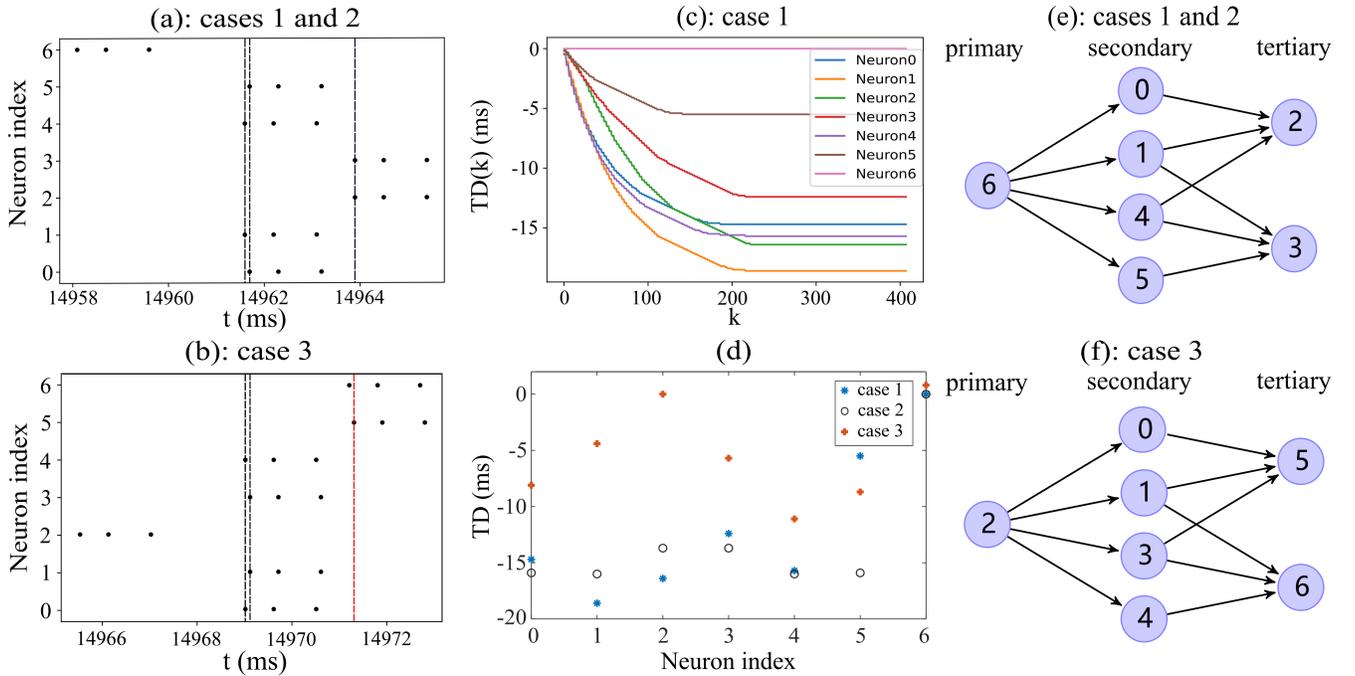}
\caption{(a-b) show locked phases for cases 1-3. Note that black (red) lines highlight that neurons spike at the same (different) time. (c) shows the time difference $TD(k)$ for each neuron in the network when comparing with neurons in the uncoupled network with the same initial membrane potentials as given in case 1 in Table~\ref{tab:ini1}. (d) stabilized time difference $TD$ for all 7 neurons in cases 1-3. Note that $TD_6=0.8~ms$ in case 3. (e-f) spiking propagation graphs with one primary neuron for cases 1-3. The network coupling structure are the same as in Fig.~\ref{fig:3delta}(a).}\label{fig:regular1}
\end{figure}

We next consider more general initial membrane potentials of neurons. With some arbitrary neuron initial potentials given in Table~\ref{tab:ini1} (as cases 1-6) for the same network structure as in Fig.~\ref{fig:3delta}(a), we find more variety of phase locking states shown in Fig.~\ref{fig:regular1}(a,b), and their corresponding stabilized synchronization parameters $\delta$ are given in Table~\ref{tab:ini1}. We remark here cases 1 \& 2 share the same locked phase and synchronization parameter $\delta$ though they have different initial potentials except for the largest one (neuron 6). In both cases 1 \& 2, neuron 6 generates action potential first and sends stimuli to other neurons connected to it, which is similar to simple examples in Sec.~\ref{subsec:ini_delta}. To further understand these locked phase and the network effect in spiking prorogation, we calculate the spiking time difference $TD(k)$. Fig.~\ref{fig:regular1}(c) shows an example of the time difference $TD(k)$ for all 7 neurons in case 1. We note that neuron 6 has $TD_6(k)=0$ for all spikes, indicating that neuron 6 spikes as its own independent spiking pattern, unaffected by the network, and $TD_i<0$ for all other neurons $i\neq 6$, suggesting spiking times of these neurons are advanced due to network connections. These suggest that neuron 6 dominate network interactions and spiking propagation. Similarly, neuron 2 in case 3 which has the largest membrane potential and $TD(k)=0$ for all spikes, dominates the spiking propagation in the regular network. %Note that in case 3, $TD_6=0.0008~ms>0$.

\begin{table}[htbp]%Í¼Æ¬Î»ÖÃ£¬htbp·Ö±ð´ú±íhere, top, bottom, page
\centering%¾ÓÖÐ
\caption{Initial membrane potential values $V_i$'s (in unit of $mV$) used in regular networks with 7 neurons where each neuron is connected to its 4 nearest neighbors, as shown in Fig.~\ref{fig:3delta}(a).}
  %\footnotesize
  %\scriptsize
\begin{tabular}{c|ccccccc|c}%ËÄ¸öc´ú±í¸Ã±íÒ»¹²ËÄÁÐ£¬ÄÚÈÝÈ«²¿¾ÓÖÐ
%\toprule%µÚÒ»µÀºáÏß
\Xhline{0.8pt}
Case & $V_0$& $V_1$& $V_2$& $V_3$ & $V_4$ & $V_5$& $V_6$&$\delta(\times10^{-3})$\\
%\midrule%µÚ¶þµÀºáÏß
\hline
1 &-63.3&-69.7& -70.0& -63.4&-64.6& -55.7&-52.0& 0.7337\\
2 & -65.0& -65.0& -65.0& -65.0& -65.0& -65.0& -52.0& 0.7337\\
3 & -64.6& -60.4& -54.9& -61.6& -70.0& -69.9& -59.8&0.7296 \\
4 &-51.0&-57.8& -65.1& --57.0&-56.9& -56.5&-51.3 &0.7081\\
5 & -70.0& -67.1& -58.1& -66.1& -57.0& -60.1& -69.8 &0.6704\\
6 & -51.0& -65.0& -65.0& -51.0& -65.0& -65.0& -65.0&0.5552 \\
%\bottomrule%µÚÈýµÀºáÏß
\Xhline{0.8pt}
\end{tabular}\label{tab:ini1}
\end{table}

Inspired by the role of neuron 6 in cases 1 \& 2 and neuron 2 in case 3, which dominate the spiking in inter burst, we grade neurons into different layers. More precisely, neurons with large initial potentials together with $\lvert TD \rvert$ sufficiently small are graded as {\em primary neurons}, such as neuron 6 in case 1. We remark here that neuron with second largest membrane potential and $TD\approx 0$ (say $\lvert TD \rvert<0.1~ms$) is also considered as primary neurons, but, the neuron with low initial potential and $TD=0$ not. We then refer neurons connected to primary neurons as {\em secondary neurons} and remaining neurons that are connected to secondary neurons as {\em tertiary neurons}, etc. Thus, based on determined primary neurons and network structure, we grade neurons into different layers.

Neurons in upper layers spike earlier than those in lower layers. Based neurons network structure and ignoring interactions of neurons within the same layer, we construct a {\em spiking propagation graph} as in Fig.~\ref{fig:regular1}(e,f) to indicate the propagation of spiking from layers to layers. For example, as seen from Fig.~\ref{fig:regular1}(e), the primary neuron 6 initiate the spiking and then propagate to neurons $\{0,1,4,5\}$ (secondary neurons) which propagate to neurons $\{2,3\}$ (tertiary neurons). This spiking propagation graph gives a coarse spiking order of neurons.

For neurons in the same layer, they may or may not spike at the same time and we study spiking order of neurons among the same layer in more details. For instance, among the secondary neurons in case 1, neurons $\{0, 5\}$ (similarly neurons $\{1,4\}$) spike at the same time as seen from Fig.~\ref{fig:regular1}(a) but slightly different to neurons $\{1,4\}$. We interpret this from the spiking propagation graph in Fig.~\ref{fig:regular1}(e) that secondary neurons $\{0,5\}$ (similarly, neurons $\{1,4\}$) share the same connected primary neuron and propagate to the same number of tertiary neurons. However, neurons $\{0,5\}$ propagate to one tertiary neuron whereas neurons $\{1,4\}$ propagate to two tertiary neurons; this leads to the different spike time among them. Tertiary neurons $\{2,3\}$ in case 1 spike at the same time whereas tertiary neurons $\{5,6\}$ in case 3 spike at different time with effectively the same stimuli from secondary neurons. This is probably due to $TD_6=0.8~ms>0$ in case 3 which indicates a short delay in spiking though it has a relatively large initial potential. These suggest that spiking hierarchy of neurons as well as locked phases of neurons can be understood largely from spiking propagation graph and spiking order of neurons among the same layer can be associated with number of neurons connected from upper layers.

\begin{figure}[H]
	%\centering \includegraphics[scale=0.25]{case456.pdf}
    \centering \includegraphics[width=\columnwidth]{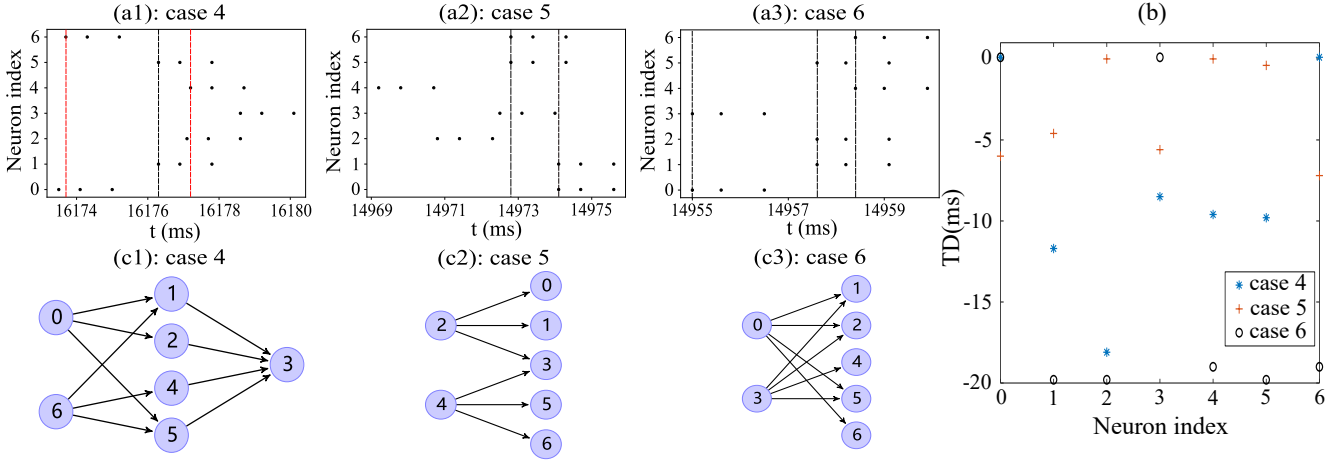}
	\caption{(a1-a3) show raster plots of locked phases; black (red) lines highlight the same (different) spiking time. (b1-b3) show the corresponding spiking propagation graphs. (c) show the stabilized spiking time different $TD$ for all neurons in cases 4-6. Network coupling structure is the same as in Fig.~\ref{fig:3delta}(a). Initial membrane potentials in each case are given in Table.~\ref{tab:ini1}.}
	\label{fig:regular2}
\end{figure}

%\subsubsection*{Neuronal network with two primary neurons}
With different initial membrane potentials, we may have two or more primary neurons. Next, we consider situations with two primary neurons in the same regular network as in Fig.~\ref{fig:3delta}(a), and see how primary neurons determine the locked phase of spiking in particular with $TD\lesssim 0$ for all neurons. More primary neurons could be studied similarly. With the regular network coupling structure in Fig.~\ref{fig:3delta}(a), there are three topologically different spiking propagation graphs with two primary neurons as shown in Fig.~\ref{fig:regular2}(c1-c3): (1) two primary neurons are connected, and all neurons are graded into three layers; (2) two primary neurons are connected, and all neurons are graded into two layers; (3) two primary neurons are not directly connected. Fig.~\ref{fig:regular2}(a1-a3) show examples of locked phases with corresponding spiking propagation graphs in Fig.~\ref{fig:regular2}(c1-c3) respectively; parameters for these locked phases are given as cases 4-6 in Table~\ref{tab:ini1}, together with Fig.~\ref{fig:regular2}(b) for the stabilized spiking time difference $TD$.

Similar to the cases with one primary neuron, neurons initiate spike in primary neurons, which propagate to secondary and then tertiary neurons. Among the same layer of neurons, they may spike at different time. In case 4, neurons $\{0,6$\} are both the primary neurons with $TD_{0,6}=0$. Neuron 0 which has a larger initial potential than neuron 6, spikes earlier than neuron 6 in inter-burst periods. Secondary neurons $\{1,5$\} spike at the same time as they are both stimulated from two primary neurons $\{0,6\}$ and propagate to the same tertiary neuron. Secondary neurons $\{2,4$\} spike at different time in inter-bursts as the only primary neuron they are connected to spikes at different time; in other words, neuron 4 receives stimulus from primary neuron 6 earlier than neuron 2 and spikes earlier than neuron 2. We remark here that if primary neurons $\{0,6\}$ have the same initial potentials, then the secondary neurons $\{2,4\}$ would have the same spiking time as well. Similarly, in case 5 secondary neurons $\{5,6\}$ spike at the same time and earlier than secondary neurons $\{0,1\}$ as their connected primary neuron $4$ spike earlier than neuron $2$; see Fig.~\ref{fig:regular2}(a2,c2). In case 6, secondary neurons $\{1,2,5\}$ spike at the same time as they receive stimuli from both primary neuron $\{0,3\}$, and spike earlier than secondary neurons $\{4,6\}$ which receive stimulus from only one primary neuron as seen in Fig.~\ref{fig:regular2}(a3,c3).

These suggest that earlier spiking time of neurons would propagate earlier and thus lead to earlier spiking of neurons connected to it in the next layer, and neurons receive stimuli from more upper-layer neurons would spike earlier than those connected to less upper-layer neurons. Moreover, we note that with the same number of neurons in the regular network, the stabilized synchronization parameters are reduced with less layers in spiking propagation graphs (see Table.~\ref{tab:ini1}); this is due to that the more clustered spiking time would lead to better synchronization.%more synchronized.

\subsection{Spiking hierarchy in small world network synchronization}\label{subsec:swn}
\begin{figure}[H]
	%\centering
	%\includegraphics[scale=0.35]{case7788.pdf}\
	\includegraphics[width=\columnwidth]{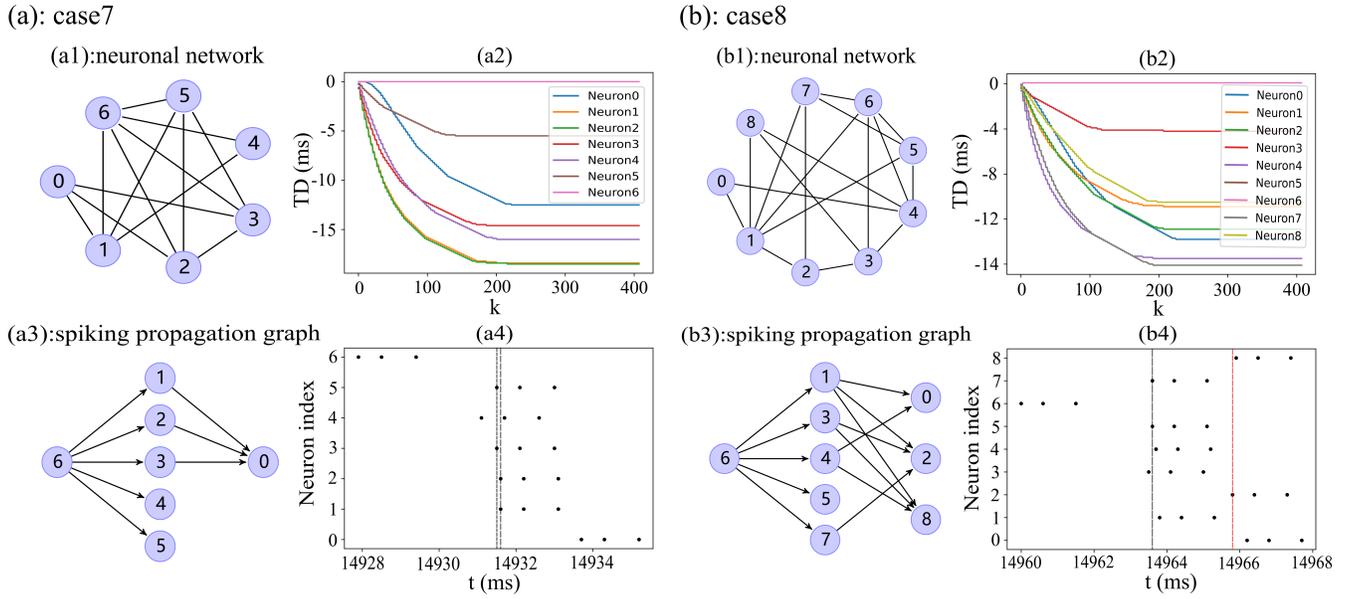}\
	\caption{(a1,b1) Topological structure of small world networks with 7 and 9 neurons as cases 7 and 8 respectively. (a2,b2) Order dependent spiking time difference $TD(k)$ compared to neurons in uncoupled network with the same initial membrane potentials given in Table~\ref{tab:ini3}. Note that the $TD(k)$ curves of neuron 5 and 7 are overlapped. (c1,c2) Spiking propagation graphs. (d1,d2) Raster plots of locked phases; black (red) lines highlight the same (different) spiking time.}
	\label{fig:random1}
\end{figure}

In this subsection, we consider small world networks which are more realistically simulated neuronal networks \cite{RS}. For simplicity, we focus on situations where $TD\lesssim 0$ for all neurons; this usually occurs when the initial potentials of primary neurons are far away from initial potentials of other neurons. We start networks with small number of neurons first. Fig.~\ref{fig:random1}(a1,b1) show two examples of small world networks with 7 and 9 neurons respectively (referred as cases 7-8). These small world networks are constructed from regular networks model with 4 nearest neighbors and each connection is rewired with probability $0.5$. Similar to cases in regular network, we construct spiking propagation graphs (shown in Fig.~\ref{fig:random1}(a3,b3)), based on neuronal initial potentials (listed in Table.~\ref{tab:ini3}), neuronal network, and their spiking time difference $TD$ compared to situations in uncoupled networks. As expected, primary neurons first spike in the inter-bursts and propagate to secondary, and then tertiary neurons.

\begin{table}[h]
\centering
\caption{Initial membrane potential values $V_i$'s (in unit of $mV$), degree $d_i's$, and received number $dr_i's$ from upper layer neurons based on the spiking propagation graph Fig.~\ref{fig:random1} for each neuron $i$.}% ±êÌâ
\begin{tabular}{c|c|c|ccccc|ccc}
%\toprule
\Xhline{0.8pt}
\multicolumn{1}{c}{} &\multicolumn{1}{c}{} & \multicolumn{1}{c}{primary}& \multicolumn{5}{c}{secondary}&\multicolumn{3}{c}{tertiary}\\
\hline
Case & $i$ &$6$&$1$&$2$&$3$&$4$&$5$&\multicolumn{3}{c}{$0$}\\
\hline
&$V_i$&-51.9&-69.7& -70.0& -63.1&-64.6& -55.7&\multicolumn{3}{c}{-63.3}  \\
7&$d_i$&5&5&5&4&2&4&\multicolumn{3}{c}{3} \\
& $dr_i$  &-&1&1&1&1&1&\multicolumn{3}{c}{3} \\
\hline
\hline
& $i$&$6$&$1$&$3$&$4$&$5$&$7$&$0$&$2$&$8$\\
\hline
& $V_i$ &-52.6&-61.6& -56.1& -68.2&-69.1& -69.1&-69.7&-67.2&-63.6 \\
8&$d_i$&5&6&4&5&4&4&2&3&3 \\
& $dr_i$  &-&1&1&1&1&1&2&3&3 \\
%\botrule
\Xhline{0.8pt}
\end{tabular}

\label{tab:ini3}
\end{table}

To examine the influence of upper-layer neurons on one neuron, we count the number of those upper-layer neurons connected to it from spiking propagation graph and denote the number as $dr$ (listed in Table.~\ref{tab:ini3}). We note the tertiary neurons $\{2,8\}$ in case 8 received stimuli from 3 secondary neurons (i.e. $dr_{2,8}=3$), which is larger than the tertiary neuron 0 received $dr_0=2$. The spiking time of these neurons can be seen from Fig.~\ref{fig:random1}(b4) which shows that neurons $\{2,8\}$ spike earlier than neuron $0$ in the inter-burst. This agrees with that more stimuli from upper layer would promote neurons to spike earlier.

In both cases 7 and 8, there are only one primary neuron, meaning secondary neurons receive stimulus from the primary neuron equally. We observe different spiking time among secondary neurons from Fig.~\ref{fig:random1}(a4,b4). In contrast to regular networks, neurons in small world networks may have different degrees. We consider the potential influence of connection degree in spiking time and list their degrees as $d$ in Table.~\ref{tab:ini3}. In case 7, secondary neurons have the degree order $d_4<d_{3,5}<d_{1,2}$. Interestingly, we note that secondary neurons spike sequentially according to their degree order in case 7. Such consistency between connection degrees and spiking hierarchy of neurons is also observed for secondary ones in case 8. This might be the minor effects of interactions in particular within the same layer stand out when interactions from the upper layer are the same. We remark here that neurons with the same degree are not necessary to spike at the same time, as slight difference among interactions within the same layer may lead to different promotion on neurons. In case 8, tertiary neurons 2 and 8 have the same degree ($d_{2,8}=3$) and receive stimuli from three secondary neurons $\{1,3,7\}$ and $\{1,3,4\}$ respectively; secondary neuron 4 spike earlier than neuron 7, which promotes neuron 2 to spike earlier than neuron 8.

%Note that neuron 3 spike earlier than neurons $\{5,7\}$ though they have the same degree in the network. This might be slight difference among interactions within secondary neurons and interactions to tertiary neurons.
%Also note that neuron 0 spike last though it has small degree than neurons $\{2,8\}$ in the tertiary neurons. This can be understood as neuron 0 receive two stimulus (i.e. $dr_0=2$) while neurons $\{2,8\}$ receive three stimulus (i.e. $dr_{2,8}=3$).

\begin{figure}[H]
	%\centering\includegraphics[scale=0.4]{201random1.pdf}
	\centering\includegraphics[width=\columnwidth]{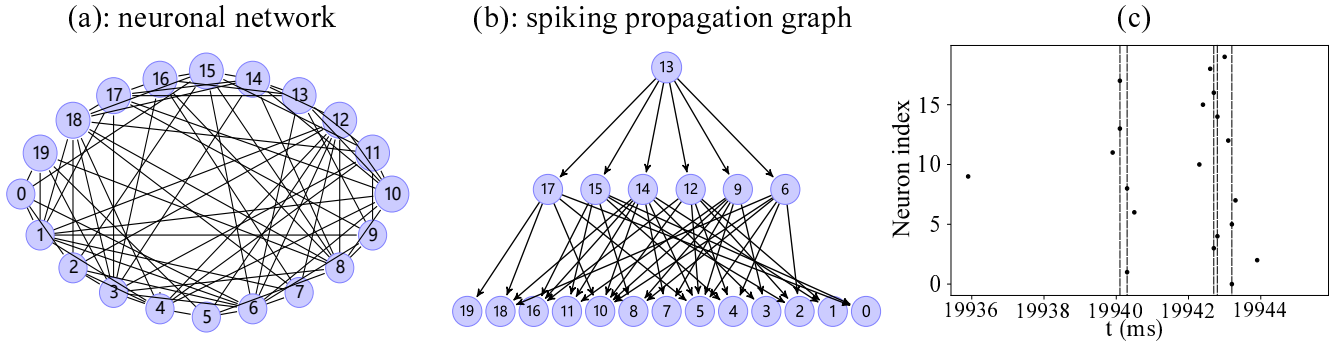}
	\caption{(a) a small world network of 20 neurons. (b) spiking propagation graph of a large small network given in (a). (c) raster plot of locked phase in this small world network; only first spiking time in the inter-burst is shown for visual effect; black (red) lines highlight the same (different) spiking time.}
	\label{fig:smallworld20}
\end{figure}

\begin{table}[h]
\centering%°Ñ±í¾ÓÖÐ
%\footnotesize
\caption{Degree ($d_i$) and spiking order of secondary neurons in a small world network given in Fig.~\ref{fig:smallworld20}(a).}% ±êÌâ
\begin{tabular}{c|cccccc}%ËÄ¸öc´ú±í¸Ã±íÒ»¹²ËÄÁÐ£¬ÄÚÈÝÈ«²¿¾ÓÖÐ
%\toprule%µÚÒ»µÀºáÏß
\Xhline{0.8pt}
 neuron $i$ &17&6&12&9&14&15 \\
%\midrule%
\hline
%$V_i$&-62.9&-62.7&-61.7&-56.4&-69.7&-63.2  \\
degree $d_i$&7&8&8&9&10&11 \\
%prediction order&1&2&3&4&4&6\\
spiking order&1&2&2&4&5&6 \\
%\bottomrule%µÚÈýµÀºáÏß
\Xhline{0.8pt}
\end{tabular}

\label{tab:ini4}
\end{table}

We next examine small world networks with large number of neurons. Fig.~\ref{fig:smallworld20} shows an example of a small world network with $20$ neurons. With some chosen initial membrane potentials, the constructed spiking propagation graph and simulated locked phase are shown in Fig.~\ref{fig:smallworld20}(b,c). Clearly, upper-layer neurons spike earlier than lower-layer neurons.
In this example, all the secondary neurons are connected to the only primary neuron. We list the degree of secondary neurons in an ascending order in Table~\ref{tab:ini4}. Clearly, neurons with smaller degree spike earlier in their inter-burst periods. Tertiary neurons are more complex than secondary neurons, as they may connected to different groups of upper-layer neurons and may also have different connection degrees. In this example, tertiary neurons have 4 different $dr$, as given in Table.~\ref{tab:ini4}. Neuron 19 receives stimuli from only one secondary neuron (with $dr_{19}=1$), and spikes latest among the tertiary neurons. Neurons $\{5,10\}$ have $dr_{5,10}=5$ and both connect with the secondary neuron $\{6,9,12,14,17\}$. Neuron 5 spikes earlier than neuron 10 which is consistent with their connection degree order $d_5<d_{10}$. Neuron $\{2,4,8,11,16,18\}$ receive stimuli from equal number of secondary neurons (though not not necessary the same group); similar for neurons $\{0,7,1,3\}$. In particular, neurons $\{2,4\}$ share the same group of connected secondary neurons, and similar to neurons $\{5,10\}$, their spiking order is consistent with their degree. Neurons $\{2,8\}$, and similarly neurons $\{11,16\}$ connected to the same number but different groups of secondary neurons. For instance, neuron 11 receives stimuli from $\{9,14,15\}$, while neuron 18 receives stimuli from $\{12,14,15\}$; earlier spiking of neuron 12 than neuron 9 together with $d_{16}<d_{11}$ leads to earlier spiking time of neuron 16.

These suggest that neurons with more stimuli (here larger $dr$) would spike earlier than others in the same layer, such as neurons 0 and 2 in case 8). For neurons in the same layer receiving effectively the same stimuli (i.e. two neurons share effectively the same spiking time of connected upper-layer neurons), those with smaller degree would then spike earlier (such as neurons 4 and 5 in case 7). For neurons in the same layer receiving the same number of stimuli and having the same degree, those receive stimuli earlier (e.g. connected upper-layer neurons spike earlier) would spike earlier, such as neurons 2 and 8 in case 8.

\begin{table}[h]
\centering%°Ñ±í¾ÓÖÐ
%\footnotesize
\caption{Degree of tertiary neurons $d_i's$, received number $dr_i's$ from upper layer neurons base on the spiking propagation graph Fig.~\ref{fig:smallworld20}(b), and spiking order for each tertiary neuron in a small world network case given in Fig.~\ref{fig:smallworld20}(a).}% ±êÌâ
\begin{tabular}{c|cc|cccccc|cccc|c}%ËÄ¸öc´ú±í¸Ã±íÒ»¹²ËÄÁÐ£¬ÄÚÈÝÈ«²¿¾ÓÖÐ
%\toprule%µÚÒ»µÀºáÏß
\Xhline{0.8pt}
neuron $i$ &$5$&$10$&$4$&$18$&$2$&$8$&$16$&$11$&$0$&$7$&$1$&$3$&$19$\\
%\midrule%
\hline
%$V_i$&\footnotesize{-64.7}&\footnotesize{-55.3}&\footnotesize{-60.0}&\footnotesize{-57.4}&
%\footnotesize{-57.7}&\footnotesize{-55.1}&\footnotesize{-62.0}&\footnotesize{-59.1}&
%\footnotesize{-69.1}&\footnotesize{-55.1}&\footnotesize{-67.0}&\footnotesize{-59.3}&\footnotesize{-54.5}\\
$dr_i$&5&5&3&3&3&3 &3&3&2&2&2&2&1\\
$d_i$&9&10&5&8&9&7&7&8&7&6&7&11&7 \\
%prediction order&1&2&3&3&5&6&6&8&9&9&9&12&13\\
spiking order&1&2&3&4&5&5&5&8&9&9&11&12&13 \\
%\bottomrule%µÚÈýµÀºáÏß
\Xhline{0.8pt}
\end{tabular}

\label{tab:ini5}
\end{table}

\section{Conclusion}\label{sec:conclusion}
In this manuscript, we study burst synchronization and the corresponding locked phase of steady spiking pattern in the aEIF neuronal network model. In particular, we focus on influence from initial membrane potentials of neurons, together with network connections. In order to measure the effect of network coupling on neuron spiking rhythm, we calculate the spiking time difference $TD$ between coupled and uncoupled networks with the same initial potentials. We focus on cases where $TD\lesssim 0$. Combined with raster plot of locked phase, we find that the neurons which have the highest initial potential always spikes first and their $TD$ approximately equal to $0$. In detail, we classify neurons in the network to primary, secondary and tertiary neurons etc. Primary neurons are grouped according to initial membrane potentials together with their spiking time difference $TD\approx 0$. Secondary neurons are those connected to primary neurons; and similarity for tertiary neurons and so on. These lead to a construction of spiking propagation graph; neurons in the upper layer spike earlier than those in lower layers. Moreover, our simulation results suggest that among the same layer, neurons receive more stimuli (larger $dr$) from upper-layer neurons would spike earlier; with effectively the same stimuli from upper-layers, neurons with smaller degree would then spike earlier; for neurons receive the same number of stimuli from upper-layer and have the same degree, if the upper-layer neurons they connected spike earlier, they would spike earlier as well.

% Our numerical simulation suggest that neurons piking order can be affected by the stimulus they received from upper-level neurons as well as their degree. Neurons receive stimulus from more up-layer neurons would spike earlier. For neurons in same layer receive stimulus from the same up-level neurons, those with less degree would spike earlier. For those with the same degree and connect to the same number of upper-level neurons, those connected with earlier spiking neurons would then spike earlier.

For large networks, situations could be complicated, neurons may receive the same number but different groups of stimuli from upper-layer, and have different degree; such as neurons $\{2, 8\}$ (similar neurons $\{1,3\}$) in the small world network case. In this case neuron 2 receive stimuli from $\{6,12,15\}$ whereas neuron 8 receive from $\{6,9,15\}$; besides the common upper-layer neurons $\{6,15\}$, neuron 2 is connected to the upper-layer neuron 12 which spikes earlier than neuron 9; however neuron 2 has larger degree than neuron 8 (i.e. $d_2>d_8$); the combined influence of degree and the upper-layer neurons in such cases remain unclear in determining their spiking hierarchy.

Our results are mainly based on $TD\lesssim 0$, and this situation occurs when initial potentials of neurons are relatively scattered, or large initial potentials are close to each other. In cases with $TD>0$ for some neurons (i.e. network connections effectively delay neurons to spike), how to determine primary neurons and then spiking propagation graph remain to be investigated in the near future. %For $TD>0$, there have

Previous studies show that increasing the number of connections in the network can enhance its synchronization \cite{BPL,CGF}. Our study provides an alternative explanation for this. With equally fixed number of neurons, networks with more connections would lead to less layers in the spiking propagation graph as more neurons are likely to connect with primary or secondary neurons. As neurons in the same layer spike at similar times, small synchronization parameter $\delta$ would be expected. In our regular network cases, we do observe less layers (e.g. cases 5 \& 6) in the spiking propagation graphs have smaller synchronization parameter (given in Table.~\ref{tab:ini1}).

%Changing the initial potential of other neurons does not affect the final phase locking state on the condition that the primary neurons and their initial voltage remain unchanged. This would lead to the same spiking propagation graph with the same primary neurons.

Last but not the least, the model we consider focus on connection given by excitatory synapses. However, there are also inhibitory synapses \cite{PBLL,PIC} that inhibit the next neuron to generate the action potential. It would be interesting to extend our model with inclusion of both excitatory and inhibitory synapses to examine the influence of initial membrane potentials on the phase locking as well as spiking hierarchy of neurons.

% Chimera states arise in a multitude of coupled systems, they have been extensively studied in a wide range of systems\cite{TRTM}. It's interesting to study the synchronization processes for chimeras like \cite{ARM}.

%In actual neural networks, there are not only excitatory synapses but also inhibitory synapses. Considering the two together in the network can improve our model. In the nervous system, the release of neurotransmitters in synapses and the random switching of ion channels all produce noise, so it is of great significance to add noise to our nonlinear system.\cite{RN}
\section*{Acknowledgments}
C. L. acknowledges financial support from National foundation of Science in China (NFSC, grant Nos. 12171179, 11871061 and 11701201) and Natural Science Foundation of Hubei Province (grant No.2020CFB847). Y. Zhang is supported by NSFC No. 11871262 and 11701200, and Hubei Key Laboratory of Engineering Modeling and Scientific Computing in HUST.

\section*{Declarations}
The authors declare no competing interests.

%\section*{Declarations}

%Some journals require declarations to be submitted in a standardised format. Please check the Instructions for Authors of the journal to which you are submitting to see if you need to complete this section. If yes, your manuscript must contain the following sections under the heading `Declarations':
%\bibliography{spikinghierarchy}
%\input{spikinghierarchy.tex}

%{\color{red}NEED TO CHECK REFERENCES. in particular, author names, abbreviation for journals}

\end{document}